\documentclass[11pt,reqno]{amsart}

\usepackage{esint}

\usepackage{cite}

\usepackage{amsmath}
\usepackage{amsthm}
\usepackage{mathrsfs}
\usepackage{graphicx}
\usepackage{mathtools}
\usepackage{amsfonts}
\usepackage{amssymb}
\usepackage{hyperref}
\usepackage[margin=1 in]{geometry}
\usepackage{enumerate}
\usepackage[normalem]{ulem}
\usepackage{tikz}
\usepackage{tikz-cd}
\usepackage{xcolor}
\usetikzlibrary{matrix,arrows,positioning,automata}
  \usepackage{cancel}
\usepackage{dsfont}

  \usepackage{scalerel,stackengine}
\stackMath
\newcommand\reallywidehat[1]{%
\savestack{\tmpbox}{\stretchto{%
  \scaleto{%
    \scalerel*[\widthof{\ensuremath{#1}}]{\kern-.6pt\bigwedge\kern-.6pt}%
    {\rule[-\textheight/2]{1ex}{\textheight}}%WIDTH-LIMITED BIG WEDGE
  }{\textheight}% 
}{0.5ex}}%
\stackon[1pt]{#1}{\tmpbox}%
}

\newcommand{\bx}{{\bm{x}}}
\newcommand{\bX}{{\boldsymbol X}}

\def\XXint#1#2#3{{\setbox0=\hbox{$#1{#2#3}{\int}$}
\vcenter{\hbox{$#2#3$}}\kern-.5\wd0}}

\newcommand{\rad}{\mathrm{rad}}

\theoremstyle{definition}

\def\smallnegint{\mathop{\int\mkern-13mu
        \raise.5ex\hbox{${\scriptscriptstyle\diagup}$}}\nolimits}

\def\bx{{\boldsymbol x}}

\def\ssetminus{\,\raise.4ex\hbox{$\scriptstyle\setminus$}\,}

\newcommand{\be}{\begin{equation}}
\newcommand{\ee}{\end{equation}}

\newcommand{\bc}{\begin{case}}
\newcommand{\ec}{\end{cases}}
\newcommand{\bs}{\begin{split}}
\newcommand{\es}{\end{split}}

\newcommand{\weak}{{\mathrm{w}}}

\newcommand{\Rd}{{\mathbb{R}^d}}

\renewcommand{\tilde}{\widetilde}

\DeclareMathOperator{\supp}{supp}

\def \be {\begin{equation}}
\def \ee {\end{equation}}

\renewcommand{\tilde}{\widetilde}

\newcommand{\trip}[1]{{\left\vert\kern-0.25ex\left\vert\kern-0.25ex\left\vert #1 
    \right\vert\kern-0.25ex\right\vert\kern-0.25ex\right\vert}}

\newcommand{\EE}{{\mathbb E}}

\newcommand{\RR}{{\mathbb R}}
\newcommand{\ZZ}{{\mathbb Z}}

\newcommand{\cadlag}{}% To make sure that \cadlag isn't already defined    
\def\cadlag/{c\`adl\`ag}

\newcommand{\eps}{\varepsilon}

\newcommand{\ind}{\mathbf{1}}

\newcommand{\oline}[1]{\overline{#1}}

\newcommand{\mcl}{\mathcal}
\newcommand{\mbb}{\mathbb}
\newcommand{\mbf}{\mathbf}

\newtheorem{lemma}{Lemma}
\newtheorem{proposition}{Proposition}
\newtheorem{theorem}{Theorem}

\theoremstyle{definition}

\numberwithin{equation}{section}
\numberwithin{lemma}{section}
\numberwithin{proposition}{section}
\numberwithin{theorem}{section}
\numberwithin{corollary}{section}
\numberwithin{definition}{section}
\numberwithin{remark}{section}

\title
[Error estimates for deterministic empirical approximations of prob. meas.]
{Error estimates for deterministic empirical approximations of probability measures}
\author{Benjamin Seeger}

\address{University of North Carolina at Chapel Hill \\ 318 Hanes Hall, CB \#3260 \\ Chapel Hill, NC 27599-3260 }
\email{bseeger@unc.edu}

\thanks{Partially supported by the National Science Foundation award DMS-2437066}


\date{\today}

\begin{document}

\maketitle

\begin{abstract}
	The question of optimally approximating an arbitrary probability measure in the Wasserstein distance by a discrete one with uniform weights is considered. Estimates are obtained for the optimal approximation distance, with an explicit rate of convergence to $0$ as the number of points tends to infinity that depends on the moment order, the parameter in the Wasserstein distance, and the  dimension. In certain low-dimensional regimes and for measures with unbounded support, the rates are improvements over those obtained through other methods, including through random sampling. Except for some critical cases, the rates are shown to be optimal.
\end{abstract}

\section{Introduction}

The purpose of this article is to explore the problem of optimally approximating a measure $\mu$ belonging to the space $\mcl P = \mcl P(\Rd)$ of probability measures on $\RR^d$ by a measure of the form
\begin{equation}\label{empirical}
	\mu^N_\bx := \frac{1}{N} \sum_{i=1}^N \delta_{x_i}
\end{equation}
for some $N \ge 1$ and $\bx = (x_1,x_2,\ldots, x_N) \in (\Rd)^N$. We are specifically interested in the quantity, for some $0 < p < \infty$,
\begin{equation}\label{error}
	e_{N;d,p}(\mu) := \inf_{\bx \in (\Rd)^N} \mcl W_p (\mu, \mu^N_x) .
\end{equation}
Here, for $0 < p < \infty$ and
\[
	\mu, \nu \in \mcl P_p = \mcl P_p(\Rd) := \left\{ \mu \in \mcl P : \mcl M_p(\mu) := \int_\Rd |x|^p \mu(dx) < \infty \right\},
\]
$\mcl W_p$ is defined by
\begin{equation}\label{Wasserstein}
	\mcl W_p(\mu,\nu) := \inf \left\{ \left( \iint_{\Rd \times \Rd} |x-y|^p \gamma(dx,dy) \right)^{1/p}: \gamma \in \Gamma(\mu,\nu) \right\},
\end{equation}
where $\Gamma(\mu,\nu) \subseteq \mcl P_p(\Rd \times \Rd)$ is the set of probability measures on $\Rd \times \Rd$ whose first and second marginals are respectively $\mu$ and $\nu$. For any $p \in (0, \infty)$, $\mcl P_p$ then becomes a metric space with distance $\mbf d_p := \mcl W_p^{p \wedge 1}$. We also define
\[
    \mcl P_\infty = \mcl P_\infty(\Rd) := \left\{ \mu \in \mcl P : \rad(\mu) := \inf\{R > 0 : \supp \mu \subseteq B_R(0) \} < \infty \right\}.
\]
%and, for $\mu, \nu \in \mcl P_{\infty}$,
%\[
%    \mcl W_\infty(\mu,\nu) := \inf_{\gamma \in \Gamma(\mu,\nu)} \gesssup_{(x,y) \in \Rd \times \Rd} |x-y|.
%\]
The quantity $e_{N;d,p}(\mu)$ will be studied for measures $\mu \in \mcl P_q$ with $0 < p < q \le \infty$.

Measures of the form $\mu^N_\bx$ for a fixed $\bx \in (\Rd)^N$ are sometimes called \textit{deterministic empirical measures} or \textit{uniform discrete measures}. The question of estimating $e_{N;d,p}$ thus falls within the study of constrained quantization. In the unconstrained setting, a probability measure is approximated by a measure with finite support, and the error
\begin{equation}\label{finite:support}
	f_{N;d,p}(\mu) := \inf \left\{ \mcl W_p(\mu, \nu) : \nu \in \mcl P, \;  \# \supp \nu \le N \right\}
\end{equation}
is analyzed. For fixed $N$ and $\mu$, this bounds $e_{N;d,p}(\mu)$ from below, since the weights in the approximating measures $\nu$ are not constrained. 

If $q > p$, then, for a fixed measure $\mu \in \mcl P_q$ that is not singular with respect to the Lebesgue measure, the quantity \eqref{finite:support} decays at the rate $N^{-1/d}$. In fact, the limit $\lim_{N \to \infty} N^{1/d} f_{N;d,p}(\mu)$ can be explicitly characterized in terms of the absolutely continuous part of $\mu$, as initiated by pioneering works of Fejes T\'oth \cite{FT} and Zador \cite{Z}. A general survey of such problems can be found in the work of Graf and Luschgy \cite{GL}.

From the lower bound $e_{N;d,p} \ge f_{N;d,p}$, it is immediate that, for any measure $\mu \in \mcl P_q$ with $q > p$ that is not singular with respect to the Lebesgue measure on $\Rd$,
\begin{equation}\label{optimal:lower:bound}
	e_{N;d,p}(\mu) \gtrsim  N^{-1/d},
\end{equation}
with a proportionality constant independent of $N$. The question of bounding $e_{N;d,p}$ from above has received recent attention due to the applications in statistics and computational mathematics \cite{BHM, DBB, DKS, MJ, GHMR1, GHMR2}. It has also become relevant in particle approximations of optimal transport and other mean field models \cite{BJ2, EN, KX, Smfg}, particularly in optimization problems for interacting agent systems such as in mean field optimal control and mean field games. The quantity $e_{N;d,p}$, sometimes with a different metric than $\mcl W_p$, then arises in the analysis of error estimates between the finite-agent system and its mean field description \cite{BCC, BEZ, CDLL, CJMS, CDJM, CDJS, DDJ}.

Upper bounds for $e_{N;d,p}(\mu)$ were obtained by Chevallier \cite{Chev} in arbitrary dimensions, with the one-dimensional case studied in great detail by Xu and Berger \cite{BX} and Bencheikh and Jourdain \cite{BJ}; the special case of Gaussian measures on a separable Hilbert space is considered in \cite{GHMR1} by Giles, Hefter, Mayer, and Ritter. In particular, it is proved in \cite{Chev} that, if $\mu$ has bounded support (i.e. $\mu \in \mcl P_\infty$), then
\begin{equation}\label{Chev:rates}
	e_{N;d,p}(\mu)
	\lesssim \rad(\mu) \times
	\begin{dcases}
		N^{-1/p} \vee N^{-1/d} & \text{if }d \ne p, \\
		N^{-1/p} \log(1+N)^{1/p} & \text{if } d = p.
	\end{dcases}
\end{equation}
In view of the lower bound \eqref{optimal:lower:bound} for nonsingular measures, the rate is optimal for $d > p$. It can be seen by example (see \cite{BX, BJ}) that the rates for $d < p$ are optimal as well, so that there is a discrepancy in the optimal rates between the unconstrained and uniform-weight settings in this low-dimensional regime.

The case of $\mu \in \mcl P_q$ with $q < \infty$ is also treated in \cite{Chev} using a truncation method, in which case, for $\mu \in \mcl P_q$ for $q > p$, the rates in \eqref{Chev:rates} are raised to the power $1 - \frac{p}{q}$. Recently, however, sharper rates were proved when $d > \frac{pq}{q-p}$ by Quattrocchi \cite{Q}, who in this regime obtained explicit characterizations of the lower and upper limits of $N^{1/d} e_{N;d,p}(\mu)$, in the same spirit as Zador's theorem. In particular, the upper limit is finite, while the lower limit is positive as soon as $\mu$ is nonsingular with respect to Lebesgue measure.

Upper bounds for $e_{N;d,p}$ can also be obtained via probabilistic methods. Namely, if the points $\bx = (x_1,x_2, \ldots, x_N)$ are chosen to be i.i.d. random variables $\bX = (X_1,X_2,\ldots, X_N)$ with distribution $\mu$, then an upper bound for $e_{N;d,p}(\mu)$ is
\begin{equation}\label{mean:convergence}
	\tilde e_{N;d,p}(\mu) := \left( \EE  \mcl W_p(\mu, \mu^N_\bX)^p \right)^{1/p}.
\end{equation}

The question of estimating $\tilde e_{N;d,p}(\mu)$ has been extensively studied, with earlier results by Ajtai, Koml\'os, and Tusn\'ady \cite{AKT} and Talagrand \cite{T} motivated by optimal matching; see also the survey by Ledoux \cite{Ledoux_survey} for a detailed discussion. Sharp (non-asymptotic) results and concentration inequalities were obtained by Fournier and Guillin \cite{FG}, adopting the multiscale approach introduced by Dereich, Scheutzow, and Schottstedt \cite{DSS}. We also mention \cite{AST}, in which Ambrosio, Stra, and Trevisan use a PDE approach; \cite{BobLed}, in which the one-dimensional setting is studied extensively by Bobkov and Ledoux; and \cite{B-LG}, in which Boissard and Le Gouic study the general setting of measures on a Polish space using the idea of the covering number in the spirit of Dudley \cite{D}. For more detailed discussions, we refer to \cite{FG} and the references therein; see also \cite{F2023}, where precise proportionality constants are derived for the rates obtained in \cite{FG}.

Among other things, it is proved in \cite{FG} that the error estimate $N^{-1/d}$ for $\tilde e_{N;d,p}(\mu)$ can be obtained for $\mu \in \mcl P_q$ in the regime where $d > 2p$ and $d > \frac{pq}{q-p}$, which therefore agrees with the same nonasymptotic bound of \cite[Theorem 1.7]{Q} in this regime. Meanwhile, if it so happens that $q < 2p$, then we may consider the regime $2p < d < \frac{pq}{q-p}$, in which case the error estimate for $\tilde e_{N;d,p}(\mu)$ (and therefore $e_{N;d,p}(\mu)$) obtained by \cite{FG} is of the order $N^{-1/p+1/q}$, which indicates that the rate $N^{-1/d + p/dq}$ obtained in \cite{Chev} for  $e_{N;d,p}(\mu)$ is suboptimal in that regime.

On the other hand, if $q > 2p$ (so that $2p > \frac{pq}{q-p}$) and $p < d < \frac{pq}{q-p}$, then \cite{FG} gives the error estimate $N^{-\frac{1}{2p}}$. 
%This has to do with variance bounds for \eqref{mean:convergence}. 
Indeed, if $\mu = \mu_\theta$ is the Bernoulli distribution with parameter $\theta \in (0,1)$, i.e. $\mu_\theta = \theta \delta_1 + (1-\theta) \delta_0$, then, for any $p > 0$,
\[
	\tilde e_{N; 1,p}(\mu_\theta) = \left( \EE \mcl W_p(\mu, \mu^N_\bX)^p \right)^{1/p} = \left( \EE \left| \frac{X_1 + X_2 + \cdots + X_N}{N} - \theta \right| \right)^{1/p} \approx N^{-\frac{1}{2p}}.
\]
This example indicates that the particular rate $N^{-\frac{1}{2p}}$ is sharp for $\tilde e_{N;1,p}$. On the other hand, observe that $e_{N;1,p}(\mu_\theta) \approx N^{-\frac{1}{p}}$, which can be achieved by setting $x_i = 1$ for $k := \lfloor N\theta \rfloor$ values of $i$ and $x_i = 0$ for the other $N - k$ values. This suggests that the upper bound $e_{N;d,p} \le \tilde e_{N; d,p}$ is not sharp in lower dimensions, where the variance bounds for the latter become nonnegligible. This discrepancy in the rates is studied in detail in $d=1$ in \cite{BX}, from which we draw the example of the Bernoulli distribution.

The main goal of the present work is to provide essentially optimal estimates for $e_{N;d,p}(\mu)$ for any dimension and for any $\mu \in \mcl P_q$ with $p < q$, including the critical cases $d = p$ and $d = \frac{qp}{q-p}$ not yet treated in the literature, as well as improving the known rates when  $p \le d \le \frac{qp}{q-p}$ and $q \ge 2p$. The proof, which relies on completely deterministic methods, is algorithmic in nature: given $d \ge 1$ and $0 < p < q \le \infty$, the placement of the points $x_1,x_2,\ldots, x_N$ can be determined up to knowledge of the measure $\mu$ at certain dyadic resolution scales depending only on $d,p,q$. This multiscale argument applies in all regimes, and is amenable to obtaining both asymptotic and nonasymptotic results.

The first main result is as follows:

\begin{theorem}\label{thm:optimal_quantize}
	Fix $0 < p < q \le \infty$. Then there exists a constant $C = C_{d,p,q} > 0$ such that, for all $N \ge 1$ and $\mu \in \mcl P_q(\RR^d)$,
	\[
		e_{N;d,p}(\mu)
		\le C \mcl M_q(\mu)^{1/q} \times
		\begin{dcases}
			N^{-1/p + 1/q} \vee N^{-1/d} & \text{if } d \ne \frac{pq}{q-p},\\
			N^{-1/d} \log(1+N)^{1/d} & \text{if } d = \frac{pq}{q-p}
		\end{dcases}
	\]
    (when $q = \infty$, the expression $\mcl M_q(\mu)^{1/q}$ above should be replaced with $\rad(\mu)$). Moreover, if $q < \infty$ and $d < \frac{qp}{q-p}$, then, for any $\mu \in \mcl P_q$, $\lim_{N \to \infty} N^{1/p - 1/q} e_{N;d,p}(\mu) = 0$.
\end{theorem}

Theorem \ref{thm:optimal_quantize} recovers the results of \cite[Theorem 3]{Chev} when $q = \infty$ (the critical case in that setting being $d = p$, where the logarithmic correction arises); the rates of \cite[Corollary 1]{Chev} when $q < \infty$ and $d < p$; and the rates of \cite{Q} when $d > \frac{pq}{q-p}$. In the regime $p \le d \le \frac{pq}{q-p}$, it provides strictly improved rates over those obtained in \cite[Corollary 1]{Chev}  (where the rate is $N^{-1/d + p/dq}$) and the rates in the same regime with $q > 2p$ implied by the method of randomly sampling \cite[Theorem 1]{FG} (where the rate is $N^{-1/2p}$). 

We obtain analogous results for measures with weak $q$-moments for $q < \infty$, i.e. measures belonging to
\[
	\mcl P_{q,\weak} = \mcl P_{q,\weak}(\Rd) := \left\{ \mu \in \mcl P : \mcl M_{q,\weak}(\mu) := \sup_{ \lambda \ge 0} \lambda^q \mu( \{ x \in \Rd : |x| \ge \lambda \} ) < \infty \right\}.
\]

\begin{theorem}\label{thm:optimal_quantize_weak}
	   Fix $0 < p < q < \infty$. Then there exists $C = C_{d,p,q} > 0$ such that, for all $N \ge 1$ and $\mu \in \mcl P_{q,\weak}$,
	\[
		e_{N;d,p}(\mu)
		\le C \mcl M_{q,\weak}(\mu)^{1/q} \times
		\begin{dcases}
			N^{-1/p + 1/q} \vee N^{-1/d} & \text{if } d \ne \frac{pq}{q-p},\\
			N^{-1/d} \log(1+N)^{1/p} & \text{if } d = \frac{pq}{q-p}.
		\end{dcases}
	\]
\end{theorem}

The focus of this paper is on the optimal rate of convergence of $e_{N;d,p}(\mu)$ to $0$ as $N \to \infty$. However, a separate but related question is to study what approximation method performs best for a given $\mu \in \mcl P$ and $N \ge 1$. For example, it is conceivable that, for certain measures with unbounded but concentrated support (e.g. with exponential tails), the truncation argument of \cite{Chev} gives a better bound than our multiscale approach. Investigating this would require a careful study of the proportionality constants in our analysis below, and a more complete answer perhaps depends on more particular information, for example, on whether or not $\mu$ is singular. 

Theorems \ref{thm:optimal_quantize} and \ref{thm:optimal_quantize_weak} are proved in Section \ref{sec:proof}, and the optimality of the rates is discussed in Section \ref{sec:optimality}. Throughout the paper, given $a,b \in \RR$ and a set $A$, we use the notation $a \lesssim_A b$ to mean that there exists a constant $C > 0$ depending only on the elements of $A$ such that $a \le Cb$. We write $a \lesssim b$ for the case $A = \{ d,p,q\}$.

\section{Multiscale analysis and proof of the error estimates}\label{sec:proof}

The key technique in our analysis, which originates from ideas in \cite{DSS} and is the basis for the arguments in \cite{FG}, is to control the Wasserstein distance by a certain multiscale quantity. For $n \ge 0$, define $Q_n = (-2^{n}, 2^{n}]^d$. For $\ell \ge 0$, let $\mcl D_\ell$ be the natural partition of $Q_0 = (-1,1]^d$ into $2^{\ell d}$ translates of the cube $(-2^{-\ell}, 2^{-\ell}]^d$. Set $B_1 = Q_0$ and, for $n \ge 2$, define $B_n = Q_{n} \backslash Q_{n-1} = (-2^n, 2^n]^d \backslash (-2^{n-1}, 2^{n-1}]^d$. Define the multiscale quantity
\begin{equation}\label{multiscale}
	\mcl L_p(\mu,\nu) := \sum_{n\ge 0} 2^{pn} \sum_{\ell \ge 0} 2^{-p \ell} \sum_{F \in \mcl D_\ell} |\mu(2^n F \cap B_n) - \nu(2^n F \cap B_n)|.
\end{equation}
Then, for all $\mu,\nu \in \mcl P_p$, we have \cite[Lemmas 5 and 6]{FG}
\begin{equation}\label{multiscale:bound}
	\mcl W_p(\mu,\nu)^p \lesssim_{d,p} \mcl L_p(\mu,\nu).
\end{equation}
%Observe that, for any $\bx = (x_1,x_2,\ldots, x_N) \in (\RR^d)^N$ and $A \subseteq \RR^d$,
%\[
%	\mu^N_\bx(A) = \frac{ \#\{ i : x_i \in A\}}{N},
%\]
%and so 
The question of minimizing \eqref{multiscale} reduces to assigning the correct number of points $x_i$ to the various sets $2^n F \cap B_n$. The next two results explain the technique for doing so.

\begin{lemma}\label{lem:decomposenumbers}
	Assume $a \ge 0$, $N \in \ZZ_{\ge 1}$, $k \in \ZZ_{\ge 0}$, and $\left| a - \frac{k}{N} \right| \le \frac{1}{N}$. Suppose that $a = \sum_{i=1}^m a_i$ for some $m \in \ZZ_{\ge 1}$ and $a_i \ge 0$. Then there exist $(k_i)_{i=1}^n \subseteq \ZZ_{\ge 0}$ such that
	\[
		\sum_{i=1}^m k_i = k \quad \text{and} \quad \left| a_i - \frac{k_i}{N} \right| \le \frac{1}{N} \quad \text{for all } i = 1,2,\ldots,m.
	\]
\end{lemma}

\begin{proof}
	Suppose first that $a = a_1 + a_2$. Assume without loss of generality that $\frac{k}{N} \le a \le \frac{k+1}{N}$;
	the case where $\frac{k-1}{N} \le a < \frac{k}{N}$ is similar. Let $i,j \in \ZZ_{\ge 0}$ be such that
	\[
		\frac{i}{N} \le a_1 < \frac{i+1}{N} \quad \text{and} \quad \frac{j}{N} \le a_2 < \frac{j+1}{N}.
	\]
    Then $\frac{i+j}{N} \le a < \frac{i+j+2}{N}$, and we split into cases depending on the value of $k$. We may have $k = i+j-1$, in which case we must have $a_1 = i/N$ and $a_2 = j/N$, and we may take $k_1 = i$ and $k_2 = j-1$. Otherwise, we either have $k = i+j$, in which case we take $k_1 = i$ and $k_2 = j$; or $k = i+j+1$, in which case we take $k_1 = i$ and $k_2 = j+1$.
	
	Now fix $m \ge 1$ and suppose that the statement of the lemma holds for any integer up to $m$. Assume $a = \sum_{i=1}^{m+1} a_i$ for some $a_i \ge 0$. Then, as was shown above, there exist $k_{m+1}, \tilde k \in \ZZ_{\ge 0}$ such that $k = \tilde k + k_{m+1}$,
	\[
		\left| \sum_{i=1}^{m} a_i - \frac{\tilde k}{N} \right| \le \frac{1}{N}, \quad \text{and} \quad \left| a_{m+1} - \frac{k_{m+1}}{N} \right| \le \frac{1}{N}.
	\]
	We conclude in view of the induction hypothesis.
\end{proof}

\begin{lemma}\label{lem:empirical_choices}
	Let $\mu$ be a nonnegative, locally-finite, additive set function on $\Rd$. Fix integers $N \ge 1$, $n \ge 0$, and $\oline \ell \ge 0$, and assume $N_n \in \ZZ_{\ge 0}$ is such that
	\[
		\left| \mu(B_n) - \frac{N_n}{N}\right| \le \frac{1}{N}.
	\]
	Then there exist nonnegative integers $k^\ell_F$ indexed by $\ell = 0,1,2, \ldots, \oline\ell$ and $F \in \mcl D_\ell$ such that $k^0_{B_1} = N_n$; for any $\ell = 1,2,\ldots, \oline\ell$ and $F' \in \mcl D_{\ell-1}$,
	\[
		\sum_{F \in \mcl D_{\ell}, \, F \subseteq F'} k^\ell_F = k^{\ell-1}_{F'} ;
	\]
	and, for all $\ell = 0,1,\ldots, \oline\ell$ and $F \in \mcl D_\ell$,
	\[
		\left|\mu(2^n F \cap B_n) - \frac{k^\ell_F}{N} \right| \le \frac{1}{N}.
	\]
\end{lemma}

\begin{proof}
	We first define $k^0_{B_1} := N_n$, and then, since $\mu(B_n) = \sum_{F \in \mcl D_1} \mu(2^nF \cap B_n)$,
	we can find the appropriate integers $(k^1_F)_{F \in \mcl D_1}$ by Lemma \ref{lem:decomposenumbers}. We then inductively apply that same Lemma in view of the fact that, for each $F' \in \mcl D_{\ell-1}$,
	\[
		\mu(2^n F' \cap B_n) = \sum_{F \in \mcl D_\ell, \; F \subseteq F'} \mu(2^n F \cap B_n).
	\]
\end{proof}

\begin{proof}[Proof of Theorem \ref{thm:optimal_quantize}]
	Fix $\mu \in \mcl P_q(\Rd)$. By a scaling argument, we may assume $\mcl M_q(\mu) \le 1$ if $q < \infty$ and $\supp \mu \subseteq (-1,1]^d$ if $q = \infty$. In particular, in the latter case, one can proceed to Step 2 below, taking $n_0 = 0$, and ignoring Step 3. Otherwise, if $q < \infty$, we proceed with Step 1.
    
    \textbf{Step 1: bounding the support of $\mu^N_\bx$.} Suppose $q < \infty$. Let $n_0 \in \ZZ_{\ge 1}$ be defined by
     \[
     	n_0 = n_0(N) := \inf\left\{ n \in \ZZ_{\ge 1} : \mu(\Rd \backslash Q_n) \le \frac{1}{N} \right\}.
     \]
	Observe that  $N \mapsto n_0(N)$ is nondecreasing, and, if $\mu$ has unbounded support, then $n_0(N) \nearrow \infty$ as $N \to \infty$. Now set $\tau_0 = 1$ and, for $n \ge 1$, $\tau_n := 2^{qn} \mu(\Rd \backslash Q_n)$. Then, for all $n \ge 0$, $\tau_n \le 1$ and $\lim_{n \to \infty} \tau_n = 0$. In particular, $2^{-q(n_0 - 1)} \tau_{n_0 - 1} \ge \frac{1}{N}$ (this is true even if $n_0 = 1$), and so, setting $\eps_N := 2^q\tau_{n_0(N) - 1}$, 
	\begin{equation}\label{exponent_of_n_0}
		2^{qn_0(N)} \le N \eps_N \quad \text{and} \quad \eps_N \le 2^q,
	\end{equation}
and, if $\mu$ has unbounded support, then $\lim_{N \to \infty} \eps_N = 0$.
     
     We impose that $x_i \in Q_{n_0}$ for all $i$; i.e., for any $n > n_0$, we have $\mu^N_{\bx}(B_n) = 0$, and so
	\begin{align}
		\sum_{n > n_0} 2^{pn}& \sum_{\ell \ge 0} 2^{-p\ell} \sum_{F \in \mcl D_\ell} |\mu(2^n F \cap B_n) - \mu^N_\bx(2^n F \cap B_n)| \notag \\
		&= \sum_{n > n_0} 2^{pn} \sum_{\ell \ge 0} 2^{-p\ell} \sum_{F \in \mcl D_\ell} \mu(2^n F \cap B_n) \notag \\
		&= \sum_{n > n_0} 2^{pn} \mu(B_n)\sum_{\ell \ge 0} 2^{-p \ell} \notag \\
		&\lesssim \left( \sum_{n > n_0} 2^{qn} \mu(B_n) \right)^{p/q} \left(\sum_{n > n_0} \mu(B_n) \right)^{1 - p/q} 
		\le N^{-1 + p/q} \kappa_N,
		\label{becauseofmoment_eps}
	\end{align}
	where $\kappa_N := \left( \sum_{n > n_0(N)} 2^{nq} \mu(B_n) \right)^{p/q}$.
    Note that $\kappa_N \lesssim \mcl M_q(\mu)^{p/q} \le 1$, and, if $\mu$ has unbounded support, then $\kappa_N \xrightarrow{N \to \infty} 0$.
	
    \textbf{Step 2: assigning the atoms.} Now, since $\mu^N_\bx(Q_{n_0}) = 1$ and
	\[
		\mu(Q_{n_0}) = 1 - \mu(\RR^d \backslash Q_{n_0}) \ge 1 -  \frac{1}{N},
	\]
	it follows from Lemma \ref{lem:decomposenumbers} that there exist nonnegative integers $(N_n)_{n = 0}^{n_0}$ such that
	\[
		\sum_{n=0}^{n_0} N_n = N \quad \text{and} 
		\quad \left| \mu(B_n) - \frac{N_n}{N} \right| 
		\le \frac{1}{N} \quad \text{for all } n = 0,1,2,\ldots, n_0.
	\]
	Fix $n = 0,1,2,\ldots, n_0$, and choose an integer $\ell_n = \ell_n(N) \ge 0$, to be specified later. Let $(k^\ell_F)$ be the nonnegative integers as specified by Lemma \ref{lem:empirical_choices} for $\oline\ell = \ell_n$. Then we can choose exactly $N_n$ of the points $x_i$ so that
\[
	\mu^N_\bx(2^n F \cap B_n) = \frac{k^\ell_F}{N} \quad \text{for all } \ell = 0,1,2,\ldots, \ell_n \quad \text{and} \quad F \in \mcl D_\ell.
\]
For any $\ell$, we have
\[
    \sum_{F \in \mcl D_\ell} |\mu(2^n F \cap B_n) - \mu^N_\bx(2^n F \cap B_n)|
    \le \mu(B_n) + \mu^N_\bx(B_n) \le 2\mu(B_n) + \frac{1}{N} =: \delta_n^N.
\]
Observe then that $\delta_n^N \lesssim 1$, and, if $q < \infty$, then $\delta_n^N  \lesssim 2^{-qn} + \frac{1}{N} \lesssim 2^{-qn}$ and
\begin{equation}\label{moment_sum}
	\sum_{n=0}^{n_0} 2^{qn} \delta_n^N \lesssim \sum_{n=0}^{n_0} 2^{qn} \mu(B_n) + \frac{1}{N} \sum_{n=0}^{n_0} 2^{qn} \lesssim \mcl M_q(\mu) + \frac{2^{qn_0}}{N} \lesssim 1.
\end{equation}
The estimates in Lemma \ref{lem:empirical_choices} then yield
\begin{align}\label{to_optimize}
	&\sum_{\ell \ge 0} 2^{-p\ell} \sum_{F \in \mcl D_\ell} |\mu(2^n F \cap B_n) - \mu^N_\bx(2^n F \cap B_n) \notag |\\
	&\le \sum_{\ell = 0}^{\ell_n} 2^{-p\ell} \sum_{F \in \mcl D_\ell} \left|\mu(2^n F \cap B_n) - \frac{k_F^\ell}{N}\right| +  \sum_{\ell = \ell_n}^{\infty} 2^{-p\ell} \delta_n^N \notag \\
	&\le \sum_{\ell=0}^{\ell_n} 2^{-(p -d)\ell} \frac{1}{N} + \sum_{\ell = \ell_n}^{\infty} 2^{-p\ell} \delta_n^N \lesssim \frac{\eps_{\ell_n}}{N} + 2^{-p\ell_n} \delta_n^N,
\end{align}
where
\[
	\eps_{\ell_n} =
	\begin{dcases}
		1, & d < p, \\
		\ell_n, & d = p, \\
		2^{(d-p)\ell_n}, & d > p.
	\end{dcases}
\]
Assume now that $q = \infty$. In that case $n_0 = 0$, and we write simply $\ell_0 = \ell_n$. We also bound $\delta_n^N \lesssim 1$. If $d < p$, we choose $\ell_0$ sufficiently large relative to $N$ that the first term on the right-hand side of \eqref{to_optimize} is dominant, and if $d \ge p$, we take $\ell_0 := \Bigg\lceil \frac{\log N}{d \log 2} \Bigg\rceil$. Inserting these choices into \eqref{to_optimize}, we then recover the results of \cite[Theorem 3]{Chev}, since then, in view of \eqref{multiscale} and \eqref{multiscale:bound},
\[
	\mcl W_p(\mu, \mu^N_\bx)^p \lesssim 
	\begin{dcases}
		N^{-1} & \text{if } d < p, \\
		\frac{\log(1+N)}{N} & \text{if } d = p, \text{ and}\\
        \frac{N^{\frac{d-p}{d}}}{N} + N^{-p/d} \lesssim N^{-p/d}  & \text{if } d > p.
    \end{dcases}
\]
For the rest of the proof, assume that $q < \infty$. If $d < p$, then, as before, we choose $\ell_n$ sufficiently large, relative to $n$ and $N$, so that the second term on the right-hand side of \eqref{to_optimize} is negligible. If $d = p$, then we fix $0 < \kappa < \frac{d^2}{q-d}$, use the fact that $\ell_n \lesssim 2^{\kappa \ell_n}$, and set
\[
    \ell_n = \inf \{ \ell \in \ZZ_{\ge 0} : 2^{(d+\kappa)\ell} \ge N \delta_n^N \}
    = \left\lceil \frac{\log (N\delta_n^N)}{(d+ \kappa) \log 2} \right \rceil \vee 0.
\]
Similarly, when $d > p$, we choose
\[
	\ell_n = \inf \{ \ell \in \ZZ_{\ge 0} : 2^{d\ell} \ge N \delta_n^N\} = \left\lceil \frac{\log (N\delta_n^N)}{d \log 2} \right \rceil \vee 0.
\]

\textbf{Step 3: collecting errors.} Proceeding with the case $q < \infty$, we set
\[
    \mcl L_n := 2^{pn} \sum_{\ell \ge 0} 2^{-p\ell} \sum_{F \in \mcl D_\ell} |\mu(2^n F \cap B_n) - \mu^N_\bx(2^n F \cap B_n)|,
\]

which, in view of \eqref{to_optimize} and the choices for $\ell_n$ in Step 2, satisfies, for $0 \le n \le n_0$,
\[
	\mcl L_n
	\lesssim 
	\begin{dcases}
        \frac{2^{pn}}{N} + 2^{-p(\ell_n - n)} \delta_n^N \lesssim
		\frac{2^{pn}}{N} & \text{if } d < p, \\
		%\frac{2^{pn}(n_0 - n + 1)}{N}
        \frac{2^{\kappa \ell_n} 2^{pn}}{N} + 2^{-p(\ell_n - n)} \delta_n^N
        \lesssim
        \frac{2^{pn}}{N} \vee \left[ N^{-\frac{d}{d+\kappa}} 2^{pn} (\delta_n^N)^{\frac{\kappa}{d + \kappa}} \right] & \text{if }d = p, \\
		\frac{2^{(d-p) \ell_n} 2^{pn}}{N} + 2^{-p(\ell_n - n)} \delta_n^N
        \lesssim
        \frac{2^{pn}}{N} \vee \left[ N^{-p/d} 2^{pn} (\delta_n^N)^{1 - \frac{p}{d}}\right]
		& \text{if }d > p.
	\end{dcases}
\]
In the cases where $d \ge p$ above, we used the fact that $\ell_n = 0$ only if $\delta_n^N \le \frac{1}{N}$.

If $d < p$, we have, in view of \eqref{exponent_of_n_0},
\[
	\sum_{n=0}^{n_0} \mcl L_n
	\lesssim \sum_{n=0}^{n_0} \frac{2^{pn}}{N} \lesssim \frac{2^{pn_0}}{N} \lesssim \eps_N^{p/q} N^{-1 + p/q}.
\]
If $d = p$, we use \eqref{moment_sum} to estimate
\begin{align*}
    \sum_{n=0}^{n_0} &N^{-\frac{d}{d+\kappa}} 2^{dn} (\delta_n^N)^{\frac{\kappa}{d + \kappa}}
    = N^{-\frac{d}{d+\kappa}}\sum_{n=0}^{n_0} (2^{qn} \delta_n^N )^{\frac{\kappa}{d+\kappa}} 2^{n\left(d - \frac{q\kappa}{d+\kappa} \right)}\\
    &\le N^{-\frac{d}{d+\kappa}} \left( \sum_{n=0}^{n_0} 2^{qn} \delta_n^N  \right)^{\frac{\kappa}{d+\kappa}}
    \left( \sum_{n=0}^{n_0} 2^{n\left(d + \kappa - \frac{q\kappa}{d} \right) } \right)^{\frac{d}{d +\kappa}} \lesssim N^{-\frac{d}{d+\kappa}} 2^{n_0\left(d - \frac{q\kappa}{d+\kappa} \right)} %\eps_N^{\frac{d}{q} - \frac{\kappa]{d+\kappa}}}},
\end{align*}
which implies, in view of \eqref{exponent_of_n_0}, that
\[
    \sum_{n=0}^{n_0} \mcl L_n 
    \lesssim N^{-1 + d/q} \left( \eps_N^{d/q}  \vee \eps_N^{\frac{d}{q} - \frac{\kappa}{d+\kappa}} \right) \lesssim N^{-1 + d/q}\eps_N^{\frac{d}{q} - \frac{\kappa}{d+\kappa}} .
\]
Similarly, if $d > p$, then
\begin{align*}
	& \sum_{n=0}^{n_0} N^{-p/d} 2^{pn} (\delta_n^N)^{1 - \frac{p}{d}}
	%= N^{-p/d} \sum_{n=0}^{n_0} 2^{-n(q-p - pq/d)} (2^{qn} \delta_n^N)^{1 - \frac{p}{d}}\\
	%&\le N^{-p/d} \left(\sum_{n=0}^{n_0} 2^{-n(qd/p - d - q)}\right)^{p/d} \left( \sum_{n=0}^{n_0} 2^{qn} \delta_n^N \right)^{1 - p/d} 
	\lesssim N^{-p/d} \left(\sum_{n=0}^{n_0} 2^{-n(qd/p - d - q)}\right)^{p/d}.
\end{align*}
The estimation of the final sum depends on whether $\frac{qd}{p} - d - q$ is positive, negative, or zero, and so
\[
	\sum_{n=0}^{n_0} \mcl L_n 
	\lesssim 
    \eps_N^{p/q} N^{-1 + p/q} \vee
	\begin{dcases}
		   N^{-p/d}  & \text{if } d > \frac{pq}{q-p}, \\
        	 N^{-p/d}n_0^{p/d} & \text{if } d = \frac{pq}{q-p},\\
		 N^{-p/d} 2^{(p-q + pq/d)n_0}  
		& \text{if } p < d < \frac{pq}{q-p}. 
	\end{dcases}
\]
In view of \eqref{exponent_of_n_0}, the right-hand-side is therefore controlled by $N^{-p/d}$ if $d > \frac{pq}{q-p}$; $N^{-p/d} \log(1+N)^{p/d}$ if $d = \frac{pq}{q-p}$; and $N^{-1 + p/q}  \eps_N^{p/q + p/d - 1} $ if $p < d < \frac{pq}{q-p}$. We conclude with the desired estimates by combining the bounds in all cases above with the bound \eqref{becauseofmoment_eps} from Step 1, using the uniform boundedness in $N$ of the quantities $\eps_N$ and $\kappa_N$, and appealing to \eqref{multiscale} and \eqref{multiscale:bound}. 

Suppose now that $q < \infty$ and $d < \frac{pq}{q-p}$. Then
\[
    N^{1/p-1/q}e_{N;d,p}(\mu)
    \lesssim
    \kappa_N^{1/p} + 
    \begin{dcases}
        \eps_N^{1/q + 1/d - 1/p} & \text{if } p < d < \frac{pq}{q-p}, \\
        \eps_N^{1/q - \kappa/d(d+\kappa)} & \text{if } d = p, \\
        \eps_N^{1/q} & \text{if } d < p,
    \end{dcases}
\]
which converges to $0$ as $N \to \infty$ if $\mu$ has unbounded support, since in that case both $\eps_N$ and $\kappa_N$ converge to $0$. If $\mu$ has bounded support, then, appealing to the case $q = \infty$,
\[
    N^{1/p-1/q} e_{N;d,p}(\mu) \lesssim \rad(\mu) \times
    \begin{dcases}
        N^{-1/q} & \text{if } d < p, \\
        N^{-1/q} \log(1+N)^{1/d} & \text{if } d = p, \\
        N^{1/p - 1/q - 1/d} & \text{if } p < d < \frac{pq}{q-p},
    \end{dcases}
\]
which, in all cases, converges to $0$ as $N \to \infty$.
\end{proof}

\begin{proof}[Proof of Theorem \ref{thm:optimal_quantize_weak}]
	The proof for measures with weak moments is very similar, and we only indicate the differences with the proof of Theorem \ref{thm:optimal_quantize}. First, assuming again that $\mcl M_{q,\weak}(\mu) \le 1$, we deduce that, for all $n$, $\mu(B_n) \le 2^q 2^{-qn}$. We then define $n_0$ to be the smallest integer satisfying $N \le 2^{-q} 2^{qn_0}$, and note that, rather than \eqref{exponent_of_n_0}, we instead have the bound $2^{qn_0} \le 2^{2q} N$. As a consequence, the estimate \eqref{becauseofmoment_eps} instead becomes
	\begin{align*}
		\sum_{n > n_0} 2^{pn}& \sum_{\ell \ge 0} 2^{-p\ell} \sum_{F \in \mcl D_\ell} |\mu(2^n F \cap B_n) - \mu^N_\bx(2^n F \cap B_n)| 
		\lesssim 2^{-(q-p)n_0} \lesssim N^{-1 +p/q}. 
	\end{align*}
	The beginning of Step 2 proceeds the same since now
	\[
		\mu(Q_{n_0}) = 1 - \mu(\RR^d \backslash Q_{n_0})\ge 1 - 2^q 2^{-qn_0} \ge 1 -  \frac{1}{N},
	\]
	and $\delta_n^N \lesssim 2^{-qn}$. In Step 3 the estimates in the regime $d \le p$, and the regime where $d > p$ and $q < \frac{dp}{d-p}$, lead to the nonasymptotic bound $N^{-1 + p/q}$, and, when $d = \frac{pq}{q-p}$, in which case the sequence $2^{qn} \delta_n^N$ is merely bounded, the estimate is $N^{-p/d} n_0 \lesssim \log(1+N) N^{-p/d}$.
	Finally, the case $d > \frac{pq}{q-p}$ proceeds similarly as in the proof of Theorem \ref{thm:optimal_quantize}.
\end{proof}

\section{On the optimality of the rates}\label{sec:optimality}

From the lower bound on $e_{N;d,p}(\mu)$ for nonsingular measures discussed in the Introduction, the optimality of the rates in Theorem \ref{thm:optimal_quantize} are evident in the regime $\frac{1}{d} < \frac{1}{p} - \frac{1}{q}$, for which the term $N^{-1/d}$ dominates $N^{-1/p + 1/q}$. Since the exponent is independent of $q$, the same is therefore also true of the bounds in Theorem \ref{thm:optimal_quantize_weak} in view of the inequalities
\[
    \mcl M_{q,\weak}(\mu) \le \mcl M_q(\mu) \lesssim_{q,q'} \mcl M_{q',\weak}(\mu)^{q/q'} \quad \text{for } q < q' < \infty \quad \text{and} \quad \mu \in \mcl P_{q', \weak}.
\]

The optimality of the rate $N^{-1/p + 1/q}$ when $\frac{1}{d} > \frac{1}{p} - \frac{1}{q}$ can be seen in various ways; note that it suffices in this case to consider measures on $\RR$ (i.e. $d = 1$), since these then embed into the higher dimensional spaces up to the critical dimension $\frac{pq}{q-p}$. First, if $\mu \in \mcl P_\infty(\RR)$ has bounded and disconnected support, then \cite[Remark 5.22]{BX}
\[
	\limsup_{N \to \infty} N^{1/p} e_{N;1,p}(\mu) > 0;
\]
note that this was already observed in the introduction in the context of the Bernoulli distribution. For the case of unbounded support, we point to \cite[Theorem 2.2]{BJ} and the discussion thereafter. In particular, if $\mu \in \mcl P(\RR)$ and if $q^* := \sup\{ q > 0 : \mcl M_q(\mu) < +\infty \}$ is finite, then
\[
	\lim_{N \to \infty} N^{\frac{1}{p} - \frac{1}{q} } e_{N;d,p}(\mu) = 
	\begin{dcases}
		0 & \text{if } q < q^*, \\
		+\infty & \text{if } q > q^*.
	\end{dcases}
\]
We may also assess the optimality of the rates in Theorem \ref{thm:optimal_quantize}, which hold uniformly over $\mu \in \mcl P_q$ with bounded moments, through studying the quantity
\begin{equation}\label{ebar}
	\oline e_{N;d,p,q} :=  \sup\left\{ e_{N;d,p}(\mu) : \mcl M_q(\mu)^{1/q} \le 1 \right\}.
\end{equation}

\begin{proposition}\label{prop:uniform_optimal}
	Assume that $q < \infty$ and $\frac{1}{d} > \frac{1}{p} - \frac{1}{q}$. Then
	\[
		\inf_{N \ge 1} N^{1/p - 1/q}\oline e_{N;d,p,q} > 0.
	\]
\end{proposition}

We consider measures $\mu \in \mcl P(\RR)$ as distributions of random variables on a sufficiently rich probability space $(\Omega, \mcl F, \mbb P)$. The key idea behind Proposition \ref{prop:uniform_optimal}, and, indeed, the analysis in the works \cite{BJ,BX}, is that, taking $(\Omega, \mcl F, \mbb P) = ([0,1], \mcl B([0,1]), \lambda)$ with $\lambda$ the Lebesgue measure, any $\mu \in \mcl P(\RR)$ can be identified with the increasing rearrangement $X: (0,1) \to \RR$, which, for $\omega \in (0,1)$, is given by $X(\omega) := \inf\{ t \in \RR : \mu( (-\infty, t]) \ge \omega \}$,
in which case, defining $\omega_i = \frac{i}{N}$, $i = 0,1,2,\ldots, N$, we have
\begin{equation}\label{charac_error}
	e_{N;1,p}(\mu) = \inf_{\bx \in \RR^N}\left( \sum_{i=1}^N \int_{\omega_{i-1}}^{\omega_i} |X(\omega) - x_i|^p d\omega \right)^{1/p}.
\end{equation}
This greatly simplifies the task of bounding the rate from below. In particular, we may use the following result to generate examples.

\begin{lemma}\label{lem:increasing:estimate}
	Assume $X: [0,1] \to \RR$ is increasing and $\dot X$ is decreasing and continuous on $[j/N,k/N]$ for some $0 \le j < k \le N$, and let $\mu = X_\sharp \lambda$ be the probability distribution of $X$ on $\RR$. Then
	\[
		e_{N;1,p}(\mu) \gtrsim_p \frac{1}{N^{1+1/p}} \left( \sum_{i=j+1}^k \dot X(\omega_i)^p \right)^{1/p}.
	\]
\end{lemma}
	
	\begin{proof}
	Let $x_i$ be the optimal point corresponding to the interval $[\omega_{i-1}, \omega_i]$ on the right-hand side of \eqref{charac_error}. Then $X(\omega_{i-1}) \le x_i \le X(\omega_i)$, and so there exists $s_i \in [\omega_{i-1}, \omega_i]$ such that $X(s_i) = x_i$. Then
	\begin{align*}
		e_{N;1,p}(\mu)^p
		&= \sum_{i=1}^N \int_{\omega_{i-1}}^{\omega_i} |X(\omega) - X(s_i)|^p d\omega\\
		&= \sum_{i=1}^N \left[  \int_{\omega_{i-1}}^{s_i} \left( \int_\omega^{s_i} \dot X(r)dr \right)^p d\omega + \int_{s_i}^{\omega_i} \left( \int_{s_i}^\omega \dot X(r)dr \right)^p d\omega \right] \\
		&\ge \sum_{i=j+1}^k \dot X(\omega_i)^p \left[  \int_{\omega_{i-1}}^{s_i} \left( s_i - \omega \right)^p d\omega + \int_{s_i}^{\omega_i} \left( \omega - s_i \right)^p d\omega \right] \\
		&= \frac{1}{p+1} \sum_{i=j+1}^k \dot X(\omega_i)^{p} \left[ (s_i - \omega_{i-1})^{p+1} + (\omega_i - s_i)^{p+1} \right] \\
		&\ge \frac{1}{2^{p+1}(p+1) N^{p+1}}\sum_{i=j+1}^k \dot X(\omega_i)^p.
	\end{align*}
\end{proof}

\begin{proof}[Proof of Proposition \ref{prop:uniform_optimal}]
	We assume without loss of generality that $d = 1$. Fix $\gamma > 1/q$, and set $\eps = \gamma q  - 1 > 0$. Define, for $N \ge 1$,
\[
	X(\omega) = X_N(\omega) = c_N \times
	\begin{dcases}
		- N^\gamma & \text{if } 0 \le \omega \le \frac{1}{N}, \\
		- \omega^{-\gamma} & \text{if } \frac{1}{N} < \omega \le 1,
	\end{dcases}
	\quad \text{where} \quad
	c_N^q =  \frac{ N^{-\eps} }{1 + \frac{1 - N^{-\eps}}{\eps} }.
\]
It is straightforward to check that $X_N$ is increasing, $\dot X_N$ is decreasing on $[1/N, 1]$, the probability distribution $\mu_N$ of $X_N$ is given by
\[
	\mu_N(dx) = \frac{1}{N} \delta_{-c_N N^\gamma}(dx) + \frac{c_N^{1/\gamma}}{\gamma} |x|^{- \frac{1+\gamma}{\gamma}} \ind_{(-c_N N^\gamma, -c_N)}(x) dx,
\]
and $\mcl M_q(\mu_N) = 1$. Then, by Lemma \ref{lem:increasing:estimate},
\begin{align*}
	e_{N;1,p}(\mu_N)^p \gtrsim_{p,\gamma} c_N^p \frac{1}{N^{1+p}} \sum_{i=2}^N \omega_i^{-p(1+\gamma)}
	\gtrsim N^{-\eps p/q} N^{-1 + p \gamma} \sum_{i=2}^N i^{-p(1+\gamma)} \gtrsim N^{-1 + p/q},
\end{align*}
and so $\inf_N N^{1/p - 1/q} \oline e_{N;1,p,q} \ge \inf_N N^{1/p - 1/q} e_{N;1,p}(\mu_N) \gtrsim 1$.
\end{proof}

Proposition \ref{prop:uniform_optimal} implies that, given $0 < p < q < \infty$ and $d < \frac{pq}{q-p}$, there cannot exist a single sequence $\alpha_N$ converging to $0$ as $N \to \infty$ such that
\[
	e_{N;d,p}(\mu) \le \mcl M_q(\mu)^{1/q} \alpha_N N^{-1/p + 1/q} \quad \text{for all } \mu \in \mcl P_q(\Rd) \text{ and } N \ge 1.
\]
It is however an interesting question to study how quickly $N^{1/p - 1/q} e_{N;d,p}(\mu)$ may converge to $0$ as $N \to \infty$ for a fixed measure $\mu$ in this regime.

On the other hand, the optimality of the rates in Theorem \ref{thm:optimal_quantize_weak} in this regime can be seen for a single measure $\mu \in \mcl P_{q,\weak}$: taking
\[
	\mu(dx) = q|x|^{-q - 1} \ind_{(-\infty, -1)}(x) dx \quad \text{on } \RR,
\]
whose increasing rearrangement is given by $X(\omega) = -\omega^{-1/q}$ for $\omega \in [0,1]$, we have $\mcl M_{q,\mathrm w}(\mu) = 1$ and, by Lemma \ref{lem:increasing:estimate},
\[
	e_{N;1,p}(\mu)^p \gtrsim \frac{1}{N^{1+p}} \sum_{i=1}^N \omega_i^{-p(1+1/q)} = N^{-1 + p/q} \sum_{i=1}^N i^{-p(1+q)} \ge N^{-1 + p/q}.
\]

We conclude by remarking that, for the critical dimension $d^* = \frac{pq}{q-p}$, it is difficult to assess whether or not the logarithmic term is necessary. In the special case $d = 1$ and $q = \infty$, for instance, it is easy to see by a direct argument that $\oline e_{N;1,p,\infty} \lesssim N^{-1/p} \vee N^{-1}$ for all $p > 0$. In general, seeking useful lower bounds for $\oline e_{N;d^*,p,q}$ when $d^* > 1$ is an intricate question.

\bibliographystyle{acm}
\bibliography{quant_rates}

\end{document}